\documentclass[10pt]{article}
\usepackage{amssymb}
\usepackage{graphicx}
\usepackage{xcolor} % A package to add color.
\usepackage{amsmath}
\usepackage{amsthm}
\usepackage{verbatim}
\usepackage{enumitem}
\setlist[enumerate]{leftmargin=1.5em}
\setlist[itemize]{leftmargin=1.5em}
\definecolor{green}{rgb}{0,0.8,0} % Redefines the color green.
\newtheorem{theorem}{Theorem}[section]

\theoremstyle{definition}

\theoremstyle{remark}
\newtheorem{remark}[theorem]{Remark}

\numberwithin{equation}{section}
\newcommand{\nrm}[1]{\Vert#1\Vert}

\newcommand{\brk}[1]{\langle#1\rangle}

\newcommand{\nnrm}[1]{{\vert\kern-0.25ex\vert\kern-0.25ex\vert #1 
    \vert\kern-0.25ex\vert\kern-0.25ex\vert}}

\newcommand{\lap}{\Delta}

\newcommand{\rd}{\partial}
\newcommand{\nb}{\nabla}

\newcommand{\alp}{\alpha}

\newcommand{\bbR}{\mathbb R}

\newcommand{\bbT}{\mathbb T}

\newcommand{\err}{\boldsymbol{\epsilon}}		% error in the linearized electron-MHD
\newcommand{\olu}{\overline{u}}
\newcommand{\nufr}{\nu^{-\frac{1}{2}}}
\begin{document}

\title{Quasi-streamwise vortices and enhanced dissipation for the incompressible 3D Navier-Stokes equations} 
\author{In-Jee Jeong\thanks{School of Mathematics, Korea Institute for Advanced Study. E-mail: ijeong@kias.re.kr} \and Tsuyoshi Yoneda\thanks{Department of Mathematics, University of Tokyo. E-mail: yoneda@ms.u-tokyo.ac.jp}}
%\footnotetext[1]{}
%\footnotetext[2]{}
\date{\today}

%\author{In-Jee Jeong}%
%\address{KIAS, Seoul, Korea 02455}%
%\email{ijeong@kias.re.kr}%
 
%\thanks{}%
%\subjclass{}%
%\keywords{}%

%\date{\today}%
%\dedicatory{}%
%\commby{}%
% ----------------------------------------------------------------

\maketitle

% ----------------------------------------------------------------

\begin{abstract}
	We consider the 3D incompressible Navier-Stokes equations under the following $2+\frac{1}{2}$-dimensional situation: vertical vortex blob (quasi-streamwise vortices) being stretched by two-dimensional shear flow.
We prove enhanced dissipation induced by such quasi-streamwise vortices. 
\end{abstract}

%\tableofcontents

\section{Introduction}
 One of the most important law in the study of developed turbulence is 
Kolmogorov's $4/5$-law, and in the derivation of this  
law, a significant  assumption is the zeroth law (see Frisch \cite{F} for example).
The zeroth law of turbulence states that, in the limit of vanishing viscosity, the rate of kinetic energy dissipation for solutions to the incompressible Navier-Stokes equations becomes nonzero.  
To formulate this law, let us recall the 3D incompressible Navier-Stokes equations on $\mathbb{T}^3:=(\mathbb{R}/\mathbb{Z})^3$: \begin{equation} \label{NS}
\left\{
\begin{aligned}
\rd_t u^\nu + u^\nu\cdot\nabla u^\nu+\nabla p = \nu \lap u^\nu + f,\\
\nabla\cdot u^\nu = 0,\\
u^\nu(t=0) = u^\nu_0,
\end{aligned}
\right.
\end{equation} where $\nu>0$ is the viscosity and $u^\nu : \bbT^3 \rightarrow \bbR^3, p : \bbT^3 \rightarrow \bbR$ denote the velocity and pressure of the fluid, respectively. Here $f : \bbT^3 \rightarrow \bbR^3$ is some external force. Assuming that the solution is sufficiently smooth, taking the dot product of the equation with $u^\nu$ and integrating over the domain $\bbT^3$ gives the energy balance: \begin{equation}\label{eq:energy-balance}
\begin{split}
\frac{d}{dt} \frac{1}{2}\nrm{u^\nu(t)}_{L^2}^2 = \int_{\bbT^3} f(t) \cdot u^\nu(t) dx - \nu \nrm{\nabla u^\nu(t)}_{L^2}^2. 
\end{split}
\end{equation} The zeroth law then postulates that, under the normalization $\nrm{u_0^\nu}_{L^2} = 1$, the mean energy dissipation rate does not vanish as $\nu \rightarrow 0^+$: \begin{equation*}
\begin{split}
\liminf_{\nu\rightarrow 0}\nu\brk{ |\nabla u^\nu |^2 }> 0,
\end{split}
\end{equation*} where $\brk{\cdot}$ usually denotes some \textit{ensemble} or long-time, space averages. 
Laboratory experiments and numerical simulations of turbulence both confirm the above zeroth law \cite{
%BV
%,E3,KIYIU
Va}.  
%In this paper, we 
% take sequences of smooth initial data $u_0^{\nu} \in C^\infty(\bbT^3)$. 
%Hence we may take $f \equiv 0$, and the energy balance is justified. 
%However, in the short-time, a trivial version of zeroth law appears, and thus we need to avoid it carefully (see \cite{JY}). 
Keeping this classical physical phenomena in mind, in this paper, 
%we employ a simpler way to avoid it.
 we compare the Navier-Stokes flow and the heat flow with the same initial data, and   show that vortex stretching is  enhancing the viscosity dissipation. This  could be a good example clarifying a mechanism of  the actual zeroth law.
Here, by the heat flow, we shall mean the unique $L^2$ solution of  \begin{equation} \label{heat}
\left\{
\begin{aligned}
\rd_t u^{\nu,heat} = \nu\lap u^{\nu,heat},\\ 
u^{\nu,heat}(t=0) = u^\nu_0
\end{aligned}
\right.
\end{equation} in $(t,x)\in \bbR_+\times \bbT^3$.

\medskip

\noindent We now state our main result. 

\begin{theorem}
For any $\alpha\in (0,3/4)$, there is a sequence of smooth initial data $\{u_0^{\nu_n}\}_{n\ge 1}$ on $\bbT^3$ with $\nu_n = 2^{-2n}$ and $\|u_0^{\nu_n}\|_{L^2}=1$ such that, for the solution $u^{\nu_n}$ to the 3D incompressible Navier-Stokes equations \eqref{NS} with zero external force, we have 
\begin{equation}\label{enhanced-dissipation}
\liminf_{ n\rightarrow \infty}  \nu_n\int_0^{\nu_n^{2\alpha/3}}\|\nabla u^{\nu_n}(t)\|_{L^2}^2dt \ge c
\end{equation} for some absolute constant $c>0$, 
whereas for the solution to the heat equation \eqref{heat}, we have 
\begin{equation}\label{heat-equation}
\nu_n\int_0^{\nu_n^{2\alpha/3}}\|\nabla u^{\nu_n,heat}(t)\|_{L^2}^2dt\le C\nu_n^{2\alp/3} (1+\nu_n^{1-2\max\{1/2,\alpha\}} ),
\end{equation}  with some absolute constant $C>0$, so that the right hand side vanishes as $n\rightarrow\infty$. 
\end{theorem}

Recalling the energy balance \eqref{eq:energy-balance}, \eqref{enhanced-dissipation} implies that at least a fixed portion of the initial $L^2$ energy for $u^{\nu_n}$ is lost in the $O(\nu_n^{2\alp/3})$-timescale. Moreover, the comparison between \eqref{enhanced-dissipation} and \eqref{heat-equation} simply tells us that viscosity alone is not enough to obtain such an enhanced dissipation.
\begin{remark}
With a different construction, the same authors considered the zeroth law under the $2+\frac{1}{2}$-dimensional Navier-Stokes flow: small-scale horizontal vortex blob being stretched by large-scale, anti-parallel pairs of vertical vortex tubes \cite{JY}. It is  inspired by the  direct numerical simulations of the 3D Navier-Stokes equations by Goto, Saito and Kawahara \cite{GSK}. They have found that sustained turbulence in a periodic box consists of a hierarchy of antiparallel pairs of vortex tubes. 
\end{remark}

\begin{remark}
	In this paper, we set a sequence of initial data tending to highly oscillating two-dimensional shear flow, adding oscillating vorticity to the vertical direction.
	Let us mention a recent numerical simulation which have inspired this construction of  the initial data. Recently, using direct numerical simulations of the 3D Navier-Stokes equations, Motoori and Goto \cite{MG} considered the generation mechanism of  turbulent boundary layer. They found that small-scale vortices in the log layer are generated predominantly by the stretching in a strain-rate field at larger scale rather than by the mean-flow stretching. On the other hand, large-scale vortices, namely, vortices as large as the height of the $99\%$ boundary layer thickness are stretched and amplified directly by the mean flow (we call ``qusi-streamwise vortices").
In particular large-scale vortices tend to align with the stretching direction of the mean flow, which is inclined at nearly $\pi/4$ from the streamwise direction.
We are inspired by the generation mechanism of this quasi-streamwise vortices. 
\end{remark}
  
\medskip

\noindent \textbf{Notations.} We write $A\lesssim B$ if there exists some constant $C>0$ independent of $\nu$ such that $A \le CB$.
Then, we say $A \approx B$ if $A \lesssim B$ and $B\lesssim A$. Finally we write $A \simeq B$ if $A/B \rightarrow 1$ as $\nu \rightarrow 0$. 
We use the notation $\brk{\cdot,\cdot}$ for the standard $L^2$ inner product:\begin{equation*}
\begin{split}
\brk{f,g} = \int_{\bbT^3} f(x)g(x)\,dx. 
\end{split}
\end{equation*} Here $f, g$ are real-valued functions. We also have $\nrm{f}_{L^2}^2=\brk{f,f}$.

\section{Proof of the main theorem}

We consider solutions of  the following form \begin{equation*}
\begin{split}
u^{\nu}(t,x_1,x_2,x_3) = u^{L,\nu}(t,x_2) + u^{S,\nu}(t,x_1,x_2) 
\end{split}
\end{equation*} where $u^{L,\nu}$ and $u^{S,\nu}$ only have nontrivial first and third components, respectively. Furthermore, the solution depends only on $x_1$ and $x_2$. It is well known that these assumptions (so-called $2+\frac{1}{2}$ dimensional flow) propagate by the Navier-Stokes equations (cf. \cite{MB}). With some abuse of notation, we shall identify $u^{L,\nu}$ and $u^{S,\nu}$ with their single nontrivial component, and drop the dependence of the solution on $\nu$. Then, it is not difficult to see that the 3D Navier-Stokes equations reduce to two \textit{scalar} equations \begin{equation*}
\begin{split}
\rd_t u^{L} = \nu \rd_{x_2}^2 u^L
\end{split}
\end{equation*} and \begin{equation*}
\begin{split}
\rd_t u^{S} + u^L\rd_{x_1} u^{S}= \nu (\rd_{x_1}^2+\rd_{x_2}^2)u^S.
\end{split}
\end{equation*}

\subsection{Choice of $u^L$}

We take $u_0^L$ to be a rescaling of \textit{smoothed triangular} wave \begin{equation}\label{eq:tri}
\begin{split}
u_0^L(x_2) = T(\nu^{-\frac{1}{2}}x_2).
\end{split}
\end{equation} To be precise, we take $T$ to be a  $C^\infty(\mathbb{T})$-smooth function satisfying 
\begin{equation*}
T(z)=z \quad\text{for}\quad z\in \left[-\frac{1}{4}, \frac{1}{4}\right],
\end{equation*}
and $T(z)=-T(1/2+z)$. It is not difficult to see that \begin{equation*}
\begin{split}
%c \le
 \nrm{ u_0^L }_{L^2}
\approx 1.
% \le C
\end{split}
\end{equation*} 
%with constants $c,C>0$ independent of $\nu$. 
We may expand $T$ in sine series: we have
 \begin{equation*}
\begin{split}
T(z) = \sum_{i\ge 1} a_i \sin(2(2i-1)\pi z),
\end{split}
\end{equation*}
where $\{a_i\}_i$ satisfies $|a_i|\lesssim_\delta |i|^{-\delta}$ for any $\delta\ge 0$.
 Here, it is assumed that $0<\nu\le1$ is given in a way that $\nu^{-\frac{1}{2}}$ is an integer. Then, $u^L(t, x_2)$ is given by the solution to the one dimensional heat equation with initial data $u_0^L$: \begin{equation*}
\begin{split}
\rd_t u^L = \nu \rd_{x_2}^2 u^L. 
\end{split}
\end{equation*} From the explicit formula for the initial data in sine series, we have that the solution is \begin{equation*}
\begin{split}
u^L(t,x_2) =  \sum_{i\ge 1} a_i \sin\left(2(2i-1)\pi \nu^{-\frac{1}{2}}x_2\right) \exp\left( -  {4\pi^2}  (2i-1)^2 t \right).
\end{split}
\end{equation*}

We shall compare it with \begin{equation*}
\begin{split}
\overline{u}^L(t,x_2) := u_0^L(x_2) \exp\left( -  {4\pi^2} t \right).
\end{split}
\end{equation*} 
Then by direct computation (using $\{a_i\}_i$ is rapidly decaying) we see that 
\begin{equation}\label{eq:u-L-comparison}
\begin{split}
\nrm{u^L(t)- \overline{u}^L(t)}_{L^\infty} \le Ct,
\end{split}
\end{equation} for $t\ge 0$ with $C>0$ independent of $\nu>0$. Indeed, \begin{equation*}
\begin{split}
&u^L(t,x_2)- \overline{u}^L(t,x_2) 
=\\ 
&\sum_{i\ge 2} a_i \sin\left(2(2i-1)\pi \nu^{-\frac{1}{2}}x_2\right)  \left( \exp\left( -  {4\pi^2}  (2i-1)^2 t \right)
  - \exp\left( -  {4\pi^2} t \right)  \right)
\end{split}
\end{equation*} so that \begin{equation*}
\begin{split}
|u^L(t,x_2)- \overline{u}^L(t,x_2) | \le C\sum_{i\ge 2} |a_i| ( (2i-1)^2 -1 )t \le Ct. 
\end{split}
\end{equation*} for any $x_2$. The advantage of $\bar{u}^L$ over $u^L$ is that $\bar{u}^L$ is exactly linear in some region of $\bbT$, which makes it easy to analyze the corresponding transport operator.

\begin{remark}\label{large-scale-heat}
	We see that 
	\begin{equation*}
	\|\nabla u^{L,heat}(t)\|_{L^2}^2\lesssim \|\nabla u_0^L\|_{L^2}^2\lesssim \nu^{-1}.
	\end{equation*}
	Thus
	\begin{equation*}
	\nu\int_0^{ {\nu^{2\alp/3}} }\|\nabla u^{L,heat}(t)\|_{L^2}^2dt \lesssim {\nu^{2\alp/3}}.
	\end{equation*}
\end{remark}

\subsection{Choice of $u^S$}

For $\alpha\in (0,3/4)$, we assume for simplicity that $\nu^{-\alp}$ is an integer (otherwise, we can simply replace $\nu^{-\alp}$ with its integer part in the following argument). Then we set
\begin{equation}\label{eq:uS0}
\begin{split}
u^S_0(x_1,x_2) = \sin(\nu^{-\alpha}x_1) \phi(\nu^{-\frac{1}{2}}x_2).
\end{split}
\end{equation}  
Here $\phi\ge 0$ is a smooth function supported in $(0, 1/4)$ with $\phi = 1$ in $(1/8,3/16)$. Note that \begin{equation*}
\begin{split}
%c\leq 
\nrm{ u_0^S }_{L^2}\approx 1.
% \le C
\end{split}
\end{equation*}
% with $c,C>0$ independent of $\nu$.
\begin{remark}\label{small-scale-heat}
We see that
\begin{equation*}
\|\nabla u^{S,heat}(t)\|_{L^2}^2\lesssim \|\nabla u^{S}_0\|_{L^2}^2\lesssim \nu^{-2\max\{1/2,\alpha\}}.
\end{equation*}
Thus
\begin{equation*}
\nu\int_0^{{\nu^{2\alp/3}}}\|\nabla u^{S,heat}(t)\|_{L^2}^2dt\lesssim \nu^{1+2\alpha/3-2\max\{1/2,\alpha\}}.
\end{equation*}
\end{remark}

We first consider the linear advection-diffusion equation by rescaled and smoothed triangular wave, possibly with an additional term $f$:  \begin{equation}  \label{eq:r}
\left\{
\begin{aligned} 
\rd_t\overline{u}^S+ \overline{u}^L \cdot \nb \overline{u}^S = \nu \lap \overline{u}^S + f,\\
\overline{u}^S(t=0) = u^S_0. 
\end{aligned}
\right.
\end{equation}  Explicitly, the transport operator is given by \begin{equation*}
\begin{split}
\overline{u}^L \cdot \nb = e^{-c_0t}T(\nu^{-\frac{1}{2}} x_2) \rd_{x_1}, \quad c_0 = 4\pi^2. 
\end{split}
\end{equation*}
Let us consider the ansatz \begin{equation}\label{eq:uS-ansatz}
\begin{split}
\overline{u}^S(t,x_1,x_2) =   e^{-a(t)} \sin
\left(\nu^{-\alpha}(x_1 - \frac{1-e^{-c_0t}}{c_0} T(\nu^{-\frac{1}{2}}x_2))
\right)  \phi(\nu^{-\frac{1}{2}}x_2) . 
\end{split}
\end{equation} Here, $a(t)$ is a continuous function of time which will be determined below. 
We shall now see that under an appropriate choice of $a(t)$, \eqref{eq:uS-ansatz} provides an \textit{approximate} solution to \eqref{eq:r} with $f = 0$. Let us compute the error. 
Taking a time derivative to \eqref{eq:uS-ansatz}, 
\begin{equation*}
\begin{split}
&\rd_t \overline{u}^S
=
 - a'(t) \olu^S \\
&- e^{-a(t)-c_0 t} \nu^{-\alpha}T(\nu^{-\frac{1}{2}}x_2) \cos
\left(\nu^{-\alpha}(x_1 - \frac{1-e^{-c_0t}}{c_0} T(\nu^{-\frac{1}{2}}x_2))
\right) \phi(\nu^{-\frac{1}{2}}x_2)
\end{split} 
\end{equation*} so that \begin{equation*}
\begin{split}
\rd_t \olu^S + \olu^L\cdot\nb \olu^S = -a'(t)\olu^S.
\end{split}
\end{equation*}
On the other hand, applying spatial derivatives to \eqref{eq:uS-ansatz},  we have 
 \begin{equation*}
\begin{split}
\rd_{x_1x_1} \olu^S = -\nu^{-2\alpha}\olu^S, 
\end{split}
\end{equation*}\begin{equation*}
\begin{split}
\rd_{x_2} \olu^S & =  e^{-a(t)} \cos\left(\nu^{-\alpha}(x_1 - \frac{1-e^{-c_0t}}{c_0} T(\nu^{-\frac{1}{2}} x_2))\right)\times\\
& \ \ \ \ \ \left( - \frac{1-e^{-c_0t}}{c_0} \nu^{-\frac{1}{2}} \nu^{-\alpha} T'(\nufr x_2) \right)  \phi(\nu^{-\frac{1}{2}}x_2) \\
&\qquad +  e^{-a(t)} \sin
\left(\nu^{-\alpha}(x_1 - \frac{1-e^{-c_0t}}{c_0} T(\nu^{-\frac{1}{2}}x_2))
\right) \nu^{-\frac{1}{2}} \phi'(\nu^{-\frac{1}{2}}x_2)
\end{split}
\end{equation*} 
and
\begin{equation*}
\begin{split}
\rd_{x_2x_2}\olu^S & = -\olu^S \times\left( - \frac{1-e^{-c_0t}}{c_0} \nu^{-\frac{1}{2}} \nu^{-\alpha} T'(\nufr x_2) \right)^2   \\
&\qquad  +  e^{-a(t)}  \cos\left(\nu^{-\alpha}(x_1 - \frac{1-e^{-c_0t}}{c_0} T(\nu^{-\frac{1}{2}} x_2))\right)\times\\
&\qquad \ \ \ \  \left( - \frac{1-e^{-c_0t}}{c_0} \nu^{-1} \nu^{-\alpha} T''(\nufr x_2) \right)  \phi(\nu^{-\frac{1}{2}}x_2)  \\
&\qquad + 2 e^{-a(t)} \cos\left(\nu^{-\alpha}(x_1 - \frac{1-e^{-c_0t}}{c_0} T(\nu^{-\frac{1}{2}} x_2))\right)\times\\
&\qquad \ \ \ \ \  \left( - \frac{1-e^{-c_0t}}{c_0} \nu^{-\frac{1}{2}}\nu^{-\alpha}  T'(\nufr x_2) \right) \nu^{-\frac{1}{2}} \phi'(\nu^{-\frac{1}{2}}x_2) \\
&\qquad + e^{-a(t)} \sin
\left(
\nu^{-\alpha}(x_1 - \frac{1-e^{-c_0t}}{c_0} T(\nu^{-\frac{1}{2}}x_2))
\right) \nu^{-1} \phi''(\nu^{-\frac{1}{2}}x_2) .
\end{split}
\end{equation*}

Recalling that \begin{equation*}
\begin{split}
T'(z) = 1\quad\text{for}\quad z\in \left[-\frac{1}{4},\frac{1}{4}\right],
%\quad T''(z) = -8 \dlt_{ \{ z = \frac{1}{4} \} } + 8 \dlt_{ \{ z = \frac{3}{4} \} }.
\end{split}
\end{equation*} we have that $T'=1$ and $T'' = 0$ on the support of $\phi$. With these observations, $\rd_{x_2x_2}\olu^S $ is simply \begin{equation*}
\begin{split}
\rd_{x_2x_2}\olu^S & = -\olu^S \frac{(1-e^{-c_0t})^2}{c_0^2} \nu^{-1-2\alpha}     \\
&\qquad - 2\nu^{-\alp-1} e^{-a(t)} \cos\left(\nu^{-\alpha}(x_1 - \frac{1-e^{-c_0t}}{c_0} T(\nu^{-\frac{1}{2}} x_2))\right) \times\\
&\qquad\qquad \frac{1-e^{-c_0t}}{c_0}  \phi'(\nu^{-\frac{1}{2}}x_2) \\
&\qquad + \nu^{-1} e^{-a(t)} \sin
\left(
\nu^{-\alpha}(x_1 - \frac{1-e^{-c_0t}}{c_0} T(\nu^{-\frac{1}{2}}x_2))
\right)  \phi''(\nu^{-\frac{1}{2}}x_2) .
\end{split}
\end{equation*} Therefore, \begin{equation*}
\begin{split}
&\lap\olu^S=(\rd_{x_1x_1}+\rd_{x_2x_2})\olu^S = -\left(\nu^{-2\alp} + \frac{(1-e^{-c_0t})^2}{c_0^2}\nu^{-1 -2\alpha}\right) \olu^S + g,
\end{split}
\end{equation*} where \begin{equation*}
\begin{split}
g & =   - 2\nu^{-\alp-1} e^{-a(t)} \cos\left(\nu^{-\alpha}(x_1 - \frac{1-e^{-c_0t}}{c_0} T(\nu^{-\frac{1}{2}} x_2))\right)\times\\
&\qquad\qquad  \frac{1-e^{-c_0t}}{c_0}  \phi'(\nu^{-\frac{1}{2}}x_2) \\
&\qquad + \nu^{-1} e^{-a(t)} \sin
\left(\nu^{-\alpha}(x_1 - \frac{1-e^{-c_0t}}{c_0} T(\nu^{-\frac{1}{2}}x_2))
\right)  \phi''(\nu^{-\frac{1}{2}}x_2) .
\end{split}
\end{equation*} Hence, defining $a(t)$ to be the solution to the ODE \begin{equation*}
\begin{split}
a'(t) & = \nu^{1-2\alp} + \frac{(1-e^{-c_0t})^2}{c_0^2}\nu^{-2\alpha} ,\\
a(0) & = 1, 
\end{split}
\end{equation*} we have that \eqref{eq:uS-ansatz} is a solution to \begin{equation}\label{eq:oluS}
\begin{split}
\rd_t\overline{u}^S+ \overline{u}^L \cdot \nb \overline{u}^S = \nu \lap \overline{u}^S - \nu g. 
\end{split}
\end{equation}  

\begin{remark}\label{small-scale-stretching} Observe that for $O(1)$-small $t>0$, we have \begin{equation*}
	\begin{split}
	a'(t) \simeq \nu^{1-2\alp} + t^2\nu^{-2\alpha}
	\end{split}
	\end{equation*} and hence \begin{equation*}
	\begin{split}
	a(t) \simeq 1 + \nu^{1-2\alp} t + \frac{t^3}{3}\nu^{-2\alpha}. 
	\end{split}
	\end{equation*} This shows that for $\olu^S$, loss of energy occurs within an $O(\nu^{\frac{2}{3}\alpha})$-timescale. Here, it is required that 
	$\nu^{1-2\alpha}t\lesssim 1$
	for $t=\nu^{\frac{2}{3}\alpha}$,
	which is equivalent to 
	$\alp<\frac{3}{4}$. 
	Indeed, when $0<\alp<\frac{3}{4}$, we have that $\exp(-a(t)) \gtrsim 1$ for $0\le t \le \nu^{2\alp/3}$, and a direct computation gives
	\begin{equation}\label{eq:enhanced}
	\nu\int_0^{{ \nu^{2\alp/3}}}\|\partial_{x_2}\bar u^S(t)\|_{L^2}^2dt \gtrsim \nu \int_0^{ {\nu^{2\alp/3} } }  t^2 \nu^{-1-2\alp} dt \approx 1.
	\end{equation}
\end{remark}

\subsection{Difference estimate} We consider the solution of \begin{equation}\label{eq:uS}
\begin{split}
\rd_t u^S + u^L\cdot\nb u^S = \nu\lap u^S 
\end{split}
\end{equation} with initial data $u^S_0(x_1,x_2) = \sin(\nu^{-\alpha}x_1)\phi(\nu^{-\frac{1}{2}}x_2) = \olu^S_0(x_1,x_2)$ and compare it with $\olu^S$ from \eqref{eq:uS-ansatz}: the difference $\err = \olu^S-u^S$ satisfies the equation \begin{equation}\label{eq:diff}
\begin{split}
\rd_t \err + \olu^L \cdot\nb\err + (u^L-\olu^L)\cdot\nb u^S = \nu  \lap\err - \nu g. 
\end{split}
\end{equation} Taking the $L^2$ inner product with $\err$, \begin{equation*}
\begin{split}
\frac{1}{2} \frac{d}{dt} \nrm{\err}_{L^2}^2 + \nu \nrm{\nb\err}_{L^2}^2 =- \brk{\nu g,\err} - \brk{ (u^L-\olu^L)\cdot\nb u^S , \err  }
\end{split}
\end{equation*} We easily estimate \begin{equation*}
\begin{split}
\left| \brk{ (u^L-\olu^L)\cdot\nb u^S , \err  }  \right|
&\lesssim \nrm{u^L-\olu^L}_{L^\infty}\nrm{\rd_{x_1} u^S}_{L^2} \nrm{\err}_{L^2}\\
& \lesssim t \nrm{\rd_{x_1}u^S_0}_{L^2}\nrm{\err}_{L^2} \\
&  \lesssim t \nu^{-\alp}\nrm{\err}_{L^2}.  
\end{split}
\end{equation*} Moreover, observing that $\nrm{\nu g}_{L^2} \lesssim 1+t\nu^{-\alpha}
$, we have \begin{equation}\label{eq:error-H1}
\begin{split}
\left|\frac{1}{2} \frac{d}{dt} \nrm{\err}_{L^2}^2 + \nu \nrm{\nb\err}_{L^2}^2 \right| \lesssim (1+t\nu^{-\alpha})
\nrm{\err}_{L^2}. 
\end{split}
\end{equation} Therefore, \begin{equation*}
\begin{split}
\frac{d}{dt}  \nrm{\err}_{L^2} \lesssim  1+t\nu^{-\alpha} ,
\end{split}
\end{equation*} and \begin{equation}\label{eq:error}
\begin{split}
\nrm{\err(t)}_{L^2} \lesssim  t+t^2\nu^{-\alpha} .
\end{split}
\end{equation}
Finally, we have after integrating \eqref{eq:error-H1} in $t\in[0,{\nu^{2\alp/3} }]$ and applying \eqref{eq:error}, 
\begin{equation}\label{error}
\begin{split}
\nu\int_0^{{ \nu^{2\alp/3} }}\|\nabla \err(t)\|_{L^2}^2dt
&\lesssim
\int_0^{{ \nu^{2\alp/3} }}\left|\brk{\nu g,\err}\right|+\left|\brk{ (u^L-\olu^L)\cdot\nb u^S , \err  }\right|dt\\
&\lesssim \int_0^{{ \nu^{2\alp/3} }}t(1+t\nu^{-\alpha})^2dt
\lesssim  \nu^{4\alp/3} + \nu^{2\alp/3}.
\end{split}
\end{equation} 

\subsection{Completion of the proof}

We are now in a position to complete the proof. Let us briefly note that strictly speaking, we only have $\nrm{u^\nu}_{L^2} \approx 1$ while the statement of the Theorem requires $\nrm{ u^\nu}_{L^2}=1$. This can be fixed with a simple rescaling of $u^\nu$. 

By \eqref{eq:enhanced} from Remark \ref{small-scale-stretching} and \eqref{error}, 
we have \begin{equation*}
\begin{split}
\nu\int_0^{ \nu^{2\alp/3} } \nrm{ \nb u^{\nu}(t) }_{L^2}^2 dt & \ge \nu  \int_0^{ \nu^{2\alp/3} } \nrm{ \nb u^{S,\nu}(t) }_{L^2}^2 dt \\
& \geq \nu  \int_0^{ \nu^{2\alp/3} } \nrm{ \partial_{x_2}u^{S,\nu}(t) }_{L^2}^2 dt \\
&\gtrsim \nu\int_0^{ \nu^{2\alp/3} } \nrm{ \partial_{x_2} \bar u^{S, \nu}(t) }_{L^2}^2 dt - \nu  \int_0^{ \nu^{2\alp/3} } \nrm{ \partial_{x_2} \err(t) }_{L^2}^2 dt\\
%& \ge \frac{1}{2}( \nrm{u_0^{S,\nu}}_{L^2}^2 -  \nrm{u^{S,\nu}(T)}_{L^2}^2  )  \\
%& \ge \frac{1}{2}( \nrm{u_0^{S,\nu}}_{L^2}^2 -  \nrm{\overline{u}^{S,\nu}(T)}_{L^2}^2 - 2\nrm{\err}_{L^2}\nrm{u_0^{S,\nu}}_{L^2}^2  ) 
&\gtrsim
1-{ \nu^{2\alp/3} }. 
\end{split}
\end{equation*} 
This gives the enhanced dissipation statement \eqref{enhanced-dissipation}. For the heat flow case, combining Remark \ref{large-scale-heat} and Remark \ref{small-scale-heat}, we immediately obtain the desired estimate.

\subsection*{Acknowledgements}

Research of TY  was partially supported by Grant-in-Aid for Scientific Research B (17H02860, 18H01136, 18H01135 and 20H01819), Japan Society for the Promotion of Science (JSPS).  IJ has been supported  by a KIAS Individual Grant MG066202 at Korea Institute for Advanced Study, the Science Fellowship of POSCO TJ Park Foundation, and the National Research Foundation of Korea grant No. 2019R1F1A1058486.

% ----------------------------------------------------------------
\bibliographystyle{amsplain}

% ----------------------------------------------------------------

\end{document}